\documentstyle{amsart}
\newcommand{\be}{\begin{equation}}
\newcommand{\ee}{\end{equation}}
\newcommand{\beq}{\begin{eqnarray}}
\newcommand{\eeq}{\end{eqnarray}}
\newcommand{\beqn}{\begin{eqnarray*}}
\newcommand{\eeqn}{\end{eqnarray*}}

\newtheorem{Theorem}{Theorem}[section]
\newtheorem{Proposition}[Theorem]{Proposition}
\newtheorem{Lemma}[Theorem]{Lemma}
\newtheorem{Corollary}[Theorem]{Corollary}

\newtheorem{Remark}[Theorem]{Remark}

\newcommand{\F}{{\Bbb F}}

\newcommand{\Fd}{\F_2}
\newcommand{\cala}{{\cal A}}
\newcommand{\call}{{\cal L}}
\newcommand{\V}{{\Bbb V}}  %{{\cal V}}
\newcommand{\glk}{GL_k}
\newcommand{\glki}{GL_{k+i}}
\newcommand{\glf}{GL_4}

\newcommand{\otiglk}{\rlap{$\,\,\otimes$}\lower 7pt\hbox{$_{\glk}$}}
\newcommand{\otiglki}{\rlap{$\,\,\otimes$}\lower 7pt\hbox{$_{\glki}$}}
\newcommand{\otiglf}{\rlap{$\,\,\otimes$}\lower 7pt\hbox{$_{\glf}$}}

\newcommand{\otigl}{\rlap{$\,\,\otimes$}\lower 7pt\hbox{$\,_{GL}$}}
\newcommand{\otigln}[1]{\rlap{$\,\,\otimes$}\lower 7pt\hbox{$\,_{GL_{#1}}$}}
\newcommand{\otiTk}{\rlap{$\otimes$}\lower 7pt\hbox{$_{T_k}$}}
\newcommand{\oticala}{\,\rlap{$\otimes$}\lower 7pt\hbox{$_{\cala}$}\,}
\newcommand{\ov}{\overline}

\input xy
\xyoption{all}
\def\eex{{\accent"5E e}\kern-.470em\raise.3ex\hbox{\char'176}}
\def\uw{u\kern-.44em\raise.82ex\hbox{ \vrule width .12em height .0ex depth .075ex \kern-0.16em \char'56}\kern-.05em}
\def\EEX{{\accent"5E E}\kern-.60em\raise.9ex\hbox{\char'176}\kern.1em}
\def\UW{U\kern-.42em\raise1.36ex\hbox{
\vrule width .13em height .0ex depth .075ex \kern-0.16em \char'56}\kern-.07em}
\def\aah{{\accent"5E a}\kern-.62em\raise.2ex\hbox{\char'22}\kern.12em}

\begin{document}
\title{On behavior of the algebraic transfer}
\author{ROBERT R. BRUNER, L\^{E} M. H\`{A} and NGUY\EEX N H. V. H\UW NG}
\maketitle
\footnotetext[1]{%The research was supported in part by Wayne State University and
The third named author was supported in part by the Vietnam
National Research Program, Grant N$^0 140 801$.}
\begin{center}
{\it Dedicated to Professor Hu\`{y}nh M\`{u}i on the occasion of his sixtieth birthday}
\end{center}

\begin{abstract}
Let $Tr_k:\Fd\otiglk PH_i(B\V_k)\to Ext_{\cala}^{k,k+i}(\Fd,\Fd) $
be the algebraic transfer, which is defined by W. Singer as an algebraic version
of the geometrical transfer $tr_k: \pi_*^S((B\V _k)_+) \to \pi_*^S(S^0)$.
It has been shown that the algebraic transfer is highly nontrivial and, more precisely,
that $Tr_k$ is an isomorphism for $k=1, 2, 3$. However, Singer showed
that $Tr_5$ is not an epimorphism.
In this paper, we prove that $Tr_4$ does not detect the non zero element
$g_s\in Ext_{\cala}^{4,12\cdot 2^s}(\Fd,\Fd)$ for every $s\geq 1$.
As a consequence, the localized $(Sq^0)^{-1}Tr_4$ given by inverting
the squaring operation $Sq^0$ is not an epimorphism.
This gives a negative answer to a prediction by Minami.
\end{abstract}

%----------------------------------
\section{Introduction and statement of results}\label{s1}
The subject of the present paper is the algebraic transfer
$$
Tr_k:\Fd\otiglk PH_i(B\V_k)\to Ext_{\cala}^{k,k+i}(\Fd,\Fd),
$$
which is defined by W. Singer as an algebraic version
of the geometrical transfer $tr_k: \pi_*^S((B\V_k)_+) \to \pi_*^S(S^0)$
to the stable homotopy groups of spheres.
Here $\V_k$ denotes a $k-$dimensional $\Fd-$vector space, $PH_*(B\V_k)$
is the primitive part consisting of all elements in $H_*(B\V_k)$, which
are annihilated by every positive-degree operation in the mod 2 Steenrod
algebra, $\cala$. Throughout the paper, the homology is taken with
coefficients in $\Fd$.

It has been proved that $Tr_k$ is an isomorphism for $k=1,2$ by Singer~\cite{Singer89}
and for $k=3$ by Boardman~\cite{Boardman}. These data together with the fact that
$Tr=\oplus_{k\geq 0} Tr_k$ is an algebra homomorphism (see \cite{Singer89}) show
that $Tr_k$ is highly nontrivial.
Therefore, the algebraic transfer is considered to be a useful tool for studying the mysterious
cohomology of the Steenrod algebra, $Ext_{\cala}^{*,*}(\Fd,\Fd)$.
In \cite{Singer89}, Singer also gave computations to
show that $Tr_4$ is an isomorphism up to a range of internal degrees. However, he proved
that $Tr_5$ is not an epimorphism.

Based on these data, we are particularly interested in the behavior of the fourth algebraic
transfer. The following theorem is the main result of this paper.

\begin{Theorem}\label{main}
$Tr_4$ does not detect the non zero element $g_s\in Ext_{\cala}^{4,12\cdot 2^s}(\Fd,\Fd)$
for every $s\geq 1$.
\end{Theorem}
\noindent
The reader is referred to May~\cite{May} for the generator $g_1$ and to Lin~\cite{Lin}
or ~\cite{LM} for the generators $g_s$.

As a consequence, we get a negative answer to a prediction by Minami~\cite{Minami}.

\begin{Corollary}\label{localized}
The localization of the fourth algebraic transfer
$$
(Sq^0)^{-1}Tr_4 : (Sq^0)^{-1}\Fd\otiglf PH_i(B\V_4) \to (Sq^0)^{-1}Ext_{\cala}^{4,4+i}(\Fd,\Fd)
$$
given by inverting $Sq^0$ is not an epimorphism.
\end{Corollary}

It is well-known (see~\cite{Liulevicius}) that there are squaring operations
$Sq^i$ ($i \geq 0$)
acting on the cohomology of the Steenrod algebra, which share most of the properties
with $Sq^i$ on the cohomology of spaces. However, $Sq^0$ is not the identity.
We refer to Section 2 for the precise meaning of the operation
$Sq^0$ on the domain of the algebraic transfer.

We next explain the idea of the proof of Theorem~\ref{main}.

Let $P_k:=H^*(B\V_k)$ be the polynomial algebra of $k$ variables, each of degree $1$.
Then, the domain of $Tr_k$, $\Fd\otiglk PH_*(B\V_k)$, is dual to $(\Fd\oticala P_k)^{GL_k}$.
In order to prove Theorem~\ref{main}, it suffices to show that
$(\Fd\oticala P_4)^{GL_4}_{12\cdot 2^s-4}=0$, for every $s\geq 1$.

Direct calculation of
$(\Fd\oticala P_4)_{12\cdot 2^s-4}$ is difficult, as $P_4$ in degree $12\cdot 2^s-4$
is a huge  $\Fd$-vector space, e.g. its dimension is 1771 for $s=1$.
To compute it, we observe that the iterated dual squaring operation
$$
(Sq^0_*)^{s}: (\Fd\oticala P_4)_{12\cdot 2^s-4} \to (\Fd\oticala P_k)_8
$$
is an isomorphism of $GL_4$-modules for any $s\geq 1$. This isomorphism is obtained by applying
repeatedly the following proposition.

\begin{Proposition}\label{even1}
Let $k$ and $r$ be positive intergers. Suppose that
each monomial $x_1^{i_1}\cdots x_k^{i_k}$ of $P_k$
in degree $2r+k$ with at least one exponent $i_t$ even is hit. Then
$$
Sq^0_*: (\Fd\oticala P_k)_{2r+k} \to (\Fd\oticala P_k)_r
$$
is an isomorphism of  $GL_k$-modules.
\end{Proposition}
\noindent
Here, as usual, we say that a polynomial $Q$ in $P_k$ is {\it hit} if it is $\cala$-decomposable.

Further, we show that $(\Fd\oticala P_4)_8$ is an $\Fd$-vector space of dimension 55. Then,
by investigating a specific basis of it, we prove that $(\Fd\oticala P_4)^{GL_4}_8 = 0$.
As a consequence, we get $(\Fd\oticala P_4)^{GL_4}_{12\cdot 2^s-4} = 0$ for every $s\geq 1$.

The reader who does not wish to follow the invariant theory
computation above may be
satisfied by the following weaker theorem, and then would not need to
read the paper's last 3 sections.

\begin{Theorem}\label{weak}
$Tr_4$ is not an isomorphism.
\end{Theorem}

This theorem is proved by observing that, on the one hand,
$$
(\Fd\oticala P_4)^{GL_4}_{20} \cong (\Fd\oticala P_4)^{GL_4}_8,
$$
and on the other hand,
$$
Ext_{\cala}^{4,4+20}(\Fd,\Fd) = \Fd\cdot g_1 \not \cong Ext_{\cala}^{4,4+8}(\Fd,\Fd)=0.
$$

\medskip
The paper is divided into six sections and organized as follows.
Section 2 starts with a recollection of the squaring operation and ends with a proof of the
isomorphism $(\Fd\oticala P_4)_{12\cdot 2^s-4} \cong (\Fd\oticala P_k)_8$.
Theorem~\ref{weak} is proved in Section 3.
We compute $(\Fd\oticala P_4)_8$ and its $GL_4$-invariants in Section 4.
We prove Theorem~\ref{main} in Section 5. Finally, in Section 6, we describe the
$GL_4$-module structure of $(\Fd\oticala P_4)_8$.

\vspace{0.2cm}
\noindent
{\it Acknowledgment:} The research was in progress during the second named author's visit to
the IHES (France) and the third named author's visit to
Wayne State University, Detroit (Michigan) in the academic year 2001-2002.

The third named author is grateful to Daniel Frohardt, David Handel, Lowell Hansen, John Klein,
Charles McGibbon, Claude Schochet and all colleagues at the Department of
Mathematics, Wayne State University for their hospitality and for the warm working atmosphere.

The authors express their hearty thanks to Tr\aah n N. Nam for helpful discussion.

%-------------------------------------
%\section{The squaring operation}
\section{A sufficient condition for the squaring operation to be an isomorphism}

This section starts with a recollection of Kameko's squaring operation
$$
Sq^0: \Fd\otiglk PH_*(B\V_k)\to \Fd\otiglk PH_*(B\V_k).
$$
The most important property of Kameko's $Sq^0$ is that it commutes with the classical $Sq^0$
on $Ext_{\cala}^*(\Fd,\Fd)$ (defined in \cite{Liulevicius}) through the algebraic transfer
(see \cite{Boardman}, \cite{Minami}).

This squaring operation is constructed as follows.

As well known, $H^*(B\V_k)$ is the polynomial algebra, $P_k := \Fd[x_1,...,x_k]$, on $k$ generators
$x_1,...,x_k$, each of degree 1. By dualizing,
\[ H_*(B\V_k)= \Gamma(a_1,\ldots, a_k)
\] is the divided power algebra
generated by $a_1,\ldots, a_k$, each of degree 1, where $a_i$ is dual to $x_i\in H^1(B\V_k)$.
Here the duality is taken with respect to the basis of $H^*(B\V_k)$ consisting of all monomials
in $x_1,\ldots, x_k$.

In \cite{Kameko} Kameko defined a homomorphism
\[ \begin{array}{cccl}
\widetilde{Sq}^0:& H_*(B\V_k)& \to & H_*(B\V_k) ,\\[.2cm]
             & a_1^{(i_1)}\cdots a_k^{(i_k)} & \mapsto &
               a_1^{(2i_1+1)}\cdots a_k^{(2i_k+1)}\,,
  \end{array}
\] where $a_1^{(i_1)}\cdots a_k^{(i_k)}$ is dual to $x_1^{i_1}\cdots x_k^{i_k}$.
%%%%%%He noted that it maps the primitive part $PH_*(B\V_k)$ to itself.
The following lemma is well known. We give a proof to make the paper self-contained.

\begin{Lemma}
$\widetilde{Sq}^0$ is a $GL_k-$homomorphism.
\end{Lemma}

\begin{pf}
We use the explanation of $\widetilde{Sq}^0$ by Crabb and Hubbuck~\cite{CH}, which does not
depend on the chosen basis of $H_*(B\V_k)$.
The element $a(\V_k) = a_1\cdots a_k$ is nothing but the image of the generator of
$\Lambda^k(\V_k)$ under the (skew) symmetrization map
\[ \Lambda^k(\V_k) \to H_k(B\V_k) = \Gamma_k(\V_k) =
   (\,\underbrace{\V_k\otimes\cdots\otimes \V_k}_{k\;
           \text{{\scriptsize times}}}\,)_{S_k},
\]
where the symmetric group $S_k$ acts on $\V_k\otimes\cdots\otimes \V_k$ by permutations of the
factors. Let $c:H_*(B\V_k) \to H_*(B\V_k)$ be the degree-halving epimorphism, which is dual to
the Frobenius monomorphism $F:H^*(B\V_k) \to H^*(B\V_k)$ defined by $F(x) = x^2$ for any $x$.
We have
$$
\widetilde{Sq}^0 (c(y)) = a(\V_k)y,
$$
for $y\in H_*(B\V_k)$.
To prove that this is well defined we need to show that if $c(y)=0$, then $a(\V_k)y =0$.
Indeed, $c(y)=0$ implies $<c(y), x> = <y, x^2>=0$ for every $x\in H^*(B\V_k)$. Here $<\cdot, \cdot>$
denotes the dual pairing between $H_*(B\V_k)$ and $H^*(B\V_k)$. So, if we write
$y=\sum a_1^{(i_1)}\cdots a_k^{(i_k)}$, then there is at least one $i_t$
which is odd in each term of
the sum. Therefore,
$$
a(\V_k)y = a_1\cdots a_k(\sum a_1^{(i_1)}\cdots a_k^{(i_k)}) =0,
$$
because $a_t a_t^{(i_t)}=0$ for any odd $i_t$. So, $\widetilde{Sq}^0$ is well defined.

As $c$ is a $GL_k$-epimorphism, the map $\widetilde{Sq}^0$ is a $GL_k$-homomorphism.

The lemma is proved.
\end{pf}

Further, it is easy to see that $c Sq_*^{2t+1}=0$, $cSq_*^{2t}=Sq_*^tc$.
So we have
$$
Sq_*^{2t+1}\widetilde{Sq}^0 =0, \,\,\, Sq_*^{2t}\widetilde{Sq}^0= \widetilde{Sq}^0 Sq_*^t.
$$
(See \cite{Hung97} for an explicit proof.)
Therefore, $\widetilde{Sq}^0$ maps $PH_*(B\V_k)$ to itself.

Kameko's $Sq^0$ is defined by
$$
Sq^0= 1\otiglk \widetilde{Sq}^0: \Fd\otiglk PH_*(B\V_k) \to \Fd\otiglk PH_*(B\V_k).
$$

The dual homomorphism $\widetilde{Sq}^0_*: P_k \to P_k$ of $\widetilde{Sq}^0$ is obviously given by
$$
 \widetilde{Sq}^0_*(x_1^{j_1}\cdots x_k^{j_k}) =
               \left\{ \begin{array}{ll}
       x_1^{\frac{j_1-1}{2}}\cdots x_k^{\frac{j_k-1}{2}}, & j_1,...,j_k \,\,\mbox{\rm odd},\\
                0 ,                                         & \mbox{\rm otherwise}.
                       \end{array} \right.
$$
Hence
$$
Ker (\widetilde{Sq}^0_*) = \ov{Even},
$$
where $\ov{Even}$ denotes the vector subspace of $P_k$ spanned by all monomials
$x_1^{i_1}\cdots x_k^{i_k}$ with at least one exponent $i_t$ even.

Let $s : P_k \to P_k$ be a left inverse of $\widetilde{Sq}^0_*$
defined as follows:
$$
  s(x_1^{i_1}\cdots x_k^{i_k}) = x_1^{2i_1+1}\cdots x_k^{2i_k+1}.
$$
It should be noted that $s$ does not commute with the doubling
map on $\cala$, that is, in general
$$
  Sq^{2t} s  \neq  s Sq^t .
$$
However, in one particular circumstance we have the following.

\begin{Lemma}
Under the hypothesis of Proposition~\ref{even1}, the map
$$
\begin{array}{rcl}
 \ov{s}: (\Fd\oticala P_k)_r & \to & (\Fd\oticala P_k)_{2r+k} \\[1ex]
                   \ov{s} [X] &  = &  [sX]
\end{array}
$$
is a well-defined linear map.
\end{Lemma}

\begin{pf}
We start with an observation that
$$
Im (Sq^{2t}s - s Sq^t) \subset \ov{Even}.
$$
We prove this by showing equivalently that
$$
\widetilde{Sq}^0_* (Sq^{2t}s - s Sq^t) = 0.
$$
Indeed,
\beqn
\widetilde{Sq}^0_*(Sq^{2t}s - s Sq^t)&=& \widetilde{Sq}^0_* Sq^{2t}s -\widetilde{Sq}^0_* s Sq^t\\
                                     &=& Sq^t \widetilde{Sq}^0_* s - \widetilde{Sq}^0_* s Sq^t \\
                                     &=& Sq^t \cdot id - id \cdot Sq^t \\
                                     &=& 0.
\eeqn
As a consequence, $s$ maps $(\cala^+P_k)_r$ to $(\cala^+P_k + \ov{Even})_{2r+k}$.
Here and in what follows, $\cala^+$ denotes the submodule of $\cala$ consisting of all positive
degree operations.  Further, by the  hypothesis of Proposition~\ref{even}, we have
$$
(\cala^+P_k + \ov{Even})_{2r+k} \subset (\cala^+P_k )_{2r+k}.
$$
Hence, $s$ maps $(\cala^+P_k)_r$ to $(\cala^+P_k)_{2r+k}$.
So the map $\ov{s}$ is well-defined. Then it is a linear map, as $s$ is.

The lemma is proved.
\end{pf}

The following proposition is also numbered as
  Proposition~\ref{even1}

\begin{Proposition}\label{even}
Let $k$ and $r$ be positive intergers. Suppose that
each monomial $x_1^{i_1}\cdots x_k^{i_k}$ of $P_k$
in degree $2r+k$ with at least one exponent $i_t$ even is hit. Then
$$
Sq^0_*: (\Fd\oticala P_k)_{2r+k} \to (\Fd\oticala P_k)_r
$$
is an isomorphism of  $GL_k$-modules.
\end{Proposition}

\begin{pf}  %{pf*}{Proof of Proposition~\ref{even}}
On the one hand, we have $Sq^0_* \ov{s} = id_{(\Fd\oticala P_k)_r}$. Indeed,
from $\widetilde{Sq}^0_* s = id_{P_k}$, it implies
$$
  Sq^0_* \ov{s}[X] = Sq^0_* [sX] = [\widetilde{Sq}^0_* s X  ] = [X],
$$
for any $X$ in degree $r$ of $P_k$.

On the other hand, we have $\ov{s}Sq^0_* = id_{(\Fd\oticala P_k)_{2r+k}}$. Indeed, by the hypothesis,
any monomial with at least one even exponent represents the $0$ class in $(\Fd\oticala P_k)_{2r+k}$,
so we need only to check on the classes of monomials with all exponents odd. We have
\beqn
 \ov{s} Sq^0_* [x_1^{2i_1+1}\cdots x_k^{2i_k+1}] &=& \ov{s}[x_1^{i_1}\cdots x_k^{i_k}]  \\
                                                 &=& [s(x_1^{i_1}\cdots x_k^{i_k})]  \\
                                                 &=& [x_1^{2i_1+1}\cdots x_k^{2i_k+1}],
\eeqn
for any $x_1^{2i_1+1}\cdots x_k^{2i_k+1}$ in degree $2r+k$ of $P_k$.

Combining the above two equalities, $Sq^0_* \ov{s} = id_{(\Fd\oticala P_k)_r}$ and
$\ov{s} Sq^0_* = id_{(\Fd\oticala P_k)_{2r+k}}$,
we see that $Sq^0_*: (\Fd\oticala P_k)_{2r+k} \to (\Fd\oticala P_k)_r$ is an isomorphism
with  inverse $\ov{s}: (\Fd\oticala P_k)_r  \to  (\Fd\oticala P_k)_{2r+k}$.

The proposition is proved.
\end{pf}  %{pf*}

The target of this section is the following.

\begin{Lemma}\label{8-string}
For every positive integer $s$,
$$
(Sq^0_*)^{s}: (\Fd\oticala P_4)_{12\cdot 2^s-4} \to (\Fd\oticala P_4)_8
$$
is an isomorphism of $GL_4$-modules.
\end{Lemma}

\begin{pf}
By using Proposition~\ref{even} repeatedly, it suffices to show that
any monomial of $P_4$ in degree $m = 12\cdot 2^s - 4$ with at least one even exponent is hit.
Since $m$ is even,  the number of even exponents in such a monomial must be either 2 or 4.
If all exponents of the monomial are even, then it is hit by $Sq^1$.
Hence we need only to consider the case of a monomial  $R$  with exactly
two even exponents  (and so exactly two odd exponents).
We can write, up to a permutation of variables, $R =  x_1x_2 Q^2$,
where $Q$ is a monomial in degree $6\cdot 2^s - 3$.

Let $\chi$ be the anti-homomorphism in the Steenrod algebra.  The so-called
$\chi$-trick, which is known to Brown and Peterson in the mid-sixties, states that
$$
u Sq^i(v) \equiv \chi(Sq^i)(u) v  \,\bmod \cala^+ M,
$$
for $u, v$ in any $\cala$-algebra $M$. (See also Wood~\cite{Wood}.)
In our case, it claims that
$$
R = x_1x_2 Sq^{6\cdot 2^s -3}(Q)
$$
is hit if and only if $\chi (Sq^{6.2^s -3})(x_1x_2)Q$ is.
We will show $\chi (Sq^{6.2^s -3})(x_1x_2)=0$ for any $s>0$.

As $\cala$ is a commutative coalgebra,  $\chi$ is a homomorphism of coalgebras
(see \cite[Proposition 8.6]{MM}). Then we have the Cartan formula
$$
\chi(Sq^n)(uv) =\sum_{i+j=n} \chi(Sq^i)(u) \chi(Sq^j)(v).
$$
Furthermore, it is shown by Brown and Peterson in \cite{BP}
that
$$
\chi(Sq^i) (x_j) = \left \{
              \begin{array}{ll} x_j^{2^p} &  \mbox{if} \,\, i =  2^p  - 1 \,\mbox{for  some $p$}, \\
                                0  &  \mbox{otherwise},
              \end{array} \right.
$$
for $x_j$ in degree 1. So, in order to prove $\chi (Sq^{6.2^s -3})(x_1x_2)=0$ we need only to show
that $6\cdot 2^s -3$ can not be written in the form
$$ 6\cdot 2^s -3 = (2^a-1) + (2^b-1)$$
with $a\geq b$. Indeed, if we have this equality, then $b=0$ as the left hand size is odd.
So $6\cdot 2^s -3 = 2^a-1$, or equivalently $3\cdot 2^s = 2^{a-1}+ 1$.
As $s>0$, the left hand side is even, hence $a-1=0$. It implies $3\cdot 2^s = 2$.
This equality has no solution $s$.

The lemma is proved.
\end{pf}

%------------------------------------
\section{The fourth algebraic transfer is not an isomorphism}

The target of this section is to prove the following theorem, which is also numbered as
Theorem~\ref{weak}.

\begin{Theorem}
$$
Tr_4 : \Fd\otiglf PH_i(B\V_4) \to Ext_{\cala}^{4,4+i}(\Fd,\Fd)
$$
is not an isomorphism.
\end{Theorem}

\begin{pf}
For any $r$, we have a commutative diagram
\begin{center}
\setlength{\unitlength}{1cm}
\begin{picture}(6,3)(0,-.3)

\put(0,2.3){\makebox(0,0){$ (\Fd\otiglf PH_i(B\V_4))_r $}}
\put(1.9,2.3){\vector(1,0){1.6}}
\put(2.7,2.5){\makebox(0,0){$_{Tr_4 }$}}
\put(5.1,2.3){\makebox(0,0){$Ext_{\cala}^{4, 4+r}(\Fd,\Fd) $}}

\put(0,1.8){\vector(0,-1){1.3}}
\put(5,1.8){\vector(0,-1){1.3}}
\put(.4,1.2){\makebox(0,0){$_{Sq^0}$}}
\put(5.4,1.2){\makebox(0,0){$_{Sq^0}$}}

\put(0,0){\makebox(0,0){$(\Fd\otiglf PH_i(B\V_4))_{2r+4} $}}
\put(1.9,0){\vector(1,0){1.6}}
\put(2.7,.2){\makebox(0,0){$_{Tr_4 }$}}
\put(5.1, 0){\makebox(0,0){$Ext_{\cala}^{4, 8+2r}(\Fd,\Fd)\, ,$}}
\end{picture}
\end{center}
where the first vertical arrow is the Kameko $Sq^0$ and the second vertical one is the classical
$Sq^0$.

The dual statement of Lemma~\ref{8-string} for $s=2$ claims that
$$
Sq^0 : (\Fd\otiglf PH_i(B\V_4))_8 \to (\Fd\otiglf PH_i(B\V_4))_{20}
$$
is an isomorphism.
On the other hand, it is known  (May~\cite{May}) that
$$
Ext_{\cala}^{4,4+8}(\Fd,\Fd)= 0 \not\cong Ext_{\cala}^{4,4+20}(\Fd,\Fd) = \Fd\cdot g_1.
$$
This implies that $Tr_4$ is not an isomorphism. The theorem is proved.
\end{pf}

\begin{Remark}{\rm This proof does not show whether $Tr_4$ fails to be a monomorphism or fails
to be an epimorphism. We will see that actually $Tr_4$ is not an epimorphism in Section~\ref{notepi}
below.   }
\end{Remark}

%------------------------------------
\section{$GL_4$-invariants of the indecomposables of $P_4$ in degree 8}

From now on, let us write $x=x_1$, $y=x_2$, $z=x_3$ and $t=x_4$ and
denote the monomial $x^ay^bz^ct^d$ by $(a, b, c, d)$ for abbreviation.

\begin{Proposition}\label{basis}
$(\Fd\oticala P_4)_8$ is an $\Fd$-vector space of dimension 55 with a basis consisting
of the classes represented by the following monomials:
\vspace{.5ex}
\begin{itemize}
\item[(A)]
$( 7, 1, 0, 0 ), ( 7, 0, 1, 0 ), ( 7, 0, 0, 1 ), ( 1, 7, 0, 0 ), ( 1, 0, 7, 0 ), ( 1, 0, 0, 7 ),$ \\
$( 0, 7, 1, 0 ), ( 0, 7, 0, 1 ), ( 0, 1, 7, 0 ), ( 0, 1, 0, 7 ), ( 0, 0, 7, 1 ), ( 0, 0, 1, 7 ),$ \\

\item[(B)]
$( 3, 3, 1, 1 ), ( 3, 1, 3, 1 ), ( 3, 1, 1, 3 ), ( 1, 3, 3, 1 ), ( 1, 3, 1, 3 ), ( 1, 1, 3, 3 ),$ \\

\item[(C)]
$( 6, 1, 1, 0 ), ( 6, 1, 0, 1 ), ( 6, 0, 1, 1 ), ( 1, 6, 1, 0 ), ( 1, 6, 0, 1 ), ( 1, 1, 6, 0 ),$ \\
$( 1, 1, 0, 6 ), ( 1, 0, 6, 1 ), (1, 0, 1, 6 ), ( 0, 6, 1, 1 ), ( 0, 1, 6, 1 ), ( 0, 1, 1, 6 ),$ \\

\item[(D)]
$( 5, 3, 0, 0 ), ( 5, 0, 3, 0 ), ( 5, 0, 0, 3 ), ( 0, 5, 3, 0 ), ( 0, 5, 0, 3 ), ( 0, 0, 5, 3 ),$  \\

\item[(E)]
$( 5, 2, 1, 0 ), ( 5, 2, 0, 1 ), ( 5, 0, 2, 1 ), ( 2, 5, 1, 0 ), ( 2, 5, 0, 1 ), ( 2, 1, 5, 0 ),$ \\
$( 2, 1, 0, 5 ), ( 2, 0, 5, 1 ), ( 2, 0, 1, 5 ), ( 0, 5, 2, 1 ), ( 0, 2, 5, 1 ), ( 0, 2, 1, 5 ),$ \\

\item[(F)]
$( 5, 1, 1, 1 ), ( 1, 5, 1, 1 ), ( 1, 1, 5, 1 ), ( 1, 1, 1, 5 ),$ \\

\item[(G)]
$( 4, 2, 1, 1 ), ( 4, 1, 2, 1 ), ( 1, 4, 2, 1 ).$\\
\end{itemize}
\end{Proposition}

The proposition is proved by combining a couple of lemmas.

\begin{Lemma}\label{generator}
$(\Fd\oticala P_4)_8$ is generated by the 55 elements listed in Proposition~\ref{basis}.
\end{Lemma}

\begin{pf}
It is easy to see that every monomial
$(a, b, c, d)$ with $a, b, c, d$ all even is hit (more precisely by $Sq^1$).

The only monomials $(a, b, c, d)$ in degree 8 with at least one of $a, b, c, d$ odd are
the following up to permutations of the variables:
\beqn
 &&(7, 1, 0, 0 ),
(3, 3, 1, 1 ),
(6, 1, 1, 0 ),
(5, 3, 0, 0 ),
(5, 2, 1, 0 ),
(5, 1, 1, 1 ),
(4, 2, 1, 1 ),  \\
 &&   (4, 3, 1, 0 ),
      (3, 3, 2, 0 ),
      (3, 2, 2, 1 ).
\eeqn

The last 3 monomials and their permutations are expressed
in terms of the first 7 monomials
and their permutations as follows:
\beqn
(4, 3, 1, 0) &=& (2, 5, 1, 0) + Sq^4(1, 2, 1, 0) + Sq^2(2, 3, 1, 0),  \\
(3, 3, 2, 0) &=& (5, 2, 1, 0) + (2, 5, 1, 0) + Sq^4(2, 1, 1, 0)+ Sq^4(1, 2, 1, 0) \\
             & &  + Sq^2(3, 2, 1, 0) + Sq^2(2, 3, 1, 0) + Sq^1(3, 3, 1, 0), \\
(3, 2, 2, 1) &=& (5, 1, 1, 1) + (4, 2, 1, 1) + (4, 1, 2, 1)  \\
             & & + Sq^2(3, 1, 1, 1) + Sq^1(4, 1, 1, 1) + Sq^1(3, 2, 1, 1) + Sq^1(3, 1, 2, 1).
\eeqn

Hence, $(\Fd\oticala P_4)_8$ is generated by the following 7 monomials
and their permutations:
$$
( 7, 1, 0, 0 ),
( 3, 3, 1, 1 ),
( 6, 1, 1, 0 ),
( 5, 3, 0, 0 ),
( 5, 2, 1, 0 ),
( 5, 1, 1, 1 ),
( 4, 2, 1, 1 ).
$$

By the family of a monomial $(a, b, c, d)$ we mean the set of all monomials
which are obtained
from $(a, b, c, d)$ by permutations of the variables.

The monomials in the 7 families above which are not
in Proposition~\ref{basis} can be expressed in terms of
the 55 elements listed there as follows.  (We give only one
expression from each symmetry class.)
\beqn
(3, 5, 0, 0) &=& (5, 3, 0, 0) + Sq^4(2, 2, 0, 0) + Sq^2(3, 3, 0, 0), \\
(5, 1, 2, 0) &=& (6, 1, 1, 0) + (5, 2, 1, 0) + Sq^1(5, 1, 1, 0),  \\
%------------
(4, 1, 1, 2) &=& (4, 2, 1, 1) + (4, 1, 2, 1) + Sq^1(4, 1, 1, 1), \\
(2, 4, 1, 1) &=& (4, 2, 1, 1) + Sq^4(1, 1, 1, 1) + Sq^2(2, 2, 1, 1), \\%this has 3 symmetried formulas
(2, 1, 1, 4) &=& (4, 2, 1, 1) + (4, 1, 2, 1)  \\
             & &  + Sq^4(1, 1, 1, 1) + Sq^2(2, 1, 1, 2) + Sq^1(4, 1, 1, 1), \\
(1, 4, 1, 2) &=& (4, 2, 1, 1) + (1, 4, 2, 1)  \\                              %New
             & & + Sq^4(1, 1, 1, 1) + Sq^2(2, 2, 1, 1) + Sq^1(1, 4, 1, 1), \\ %new
(1, 2, 1, 4) &=& (4, 2, 1, 1) + (1, 4, 2, 1)  \\
             & &  + Sq^2(2, 2, 1, 1) + Sq^2(1, 2, 1, 2) + Sq^1(1, 4, 1, 1), \\
(1, 1, 4, 2) &=& (4, 1, 2, 1) + (1, 4, 2, 1)  \\                               %New
             & &  + Sq^2(2, 1, 2, 1) + Sq^2(1, 2, 2, 1) + Sq^1(1, 1, 4, 1), \\ %new
(1, 1, 2, 4) &=& (4, 1, 2, 1) + (1, 4, 2, 1) + Sq^4(1, 1, 1, 1) \\
             & &  + Sq^2(2, 1, 2, 1) + Sq^2(1, 2, 2, 1) + Sq^2(1, 1, 2, 2)+ Sq^1(1, 1, 4, 1).
\eeqn

The lemma is proved.
\end{pf}

\begin{Lemma}\label{independent}
The 55 elements listed in Proposition~\ref{basis} are linearly independent in $(\Fd\oticala P_4)_8$.
\end{Lemma}

\begin{pf}
We will use an equivalence relation defined by saying that, for two polynomials $P$ and $Q$, $P$ is
equivalent to $Q$, denoted by $P \sim Q$, if $P-Q$ is hit.

If $X$ is one of the letters from $A$ to $G$, let $X_i$ be the $i$-th element in family $X$
according to the order listed in Proposition~\ref{basis}. (This is the lexicographical
order in each family.)

Suppose there is a linear relation between the 55 elements listed there
$$
\sum_{i=1}^{12} a_i A_i + \sum_{i=1}^6 b_i B_i + \sum_{i=1}^{12} c_i C_i + \sum_{i=1}^6 d_i D_i +
\sum_{i=1}^{12} e_i E_i + \sum_{i=1}^4 f_i F_i + \sum_{i=1}^3 g_i G_i = 0,
$$
where $a_i, b_i, c_i, d_i, e_i, f_i, g_i \in \Fd$. We need to show all these coefficients are zero.
The proof is divided into 4 steps.

\vspace{0.2cm}
\noindent
{\it Step 1.}  We call a monomial {\it a spike} if each of
its exponents is of the form $2^n-1$ for some $n$. It is well
known that spikes do not appear in the expression of $Sq^i Y$
for any $i$ positive and any monomial $Y$, since the powers
$x^{2^{n}-1}$ are not hit in the one variable case. Hence, the coefficient of any spike
is zero in every linear relation in $\Fd\oticala P_k$.

Among the 55 elements of Proposition~\ref{basis}, the classes
of families $A$ and $B$ are spikes.  So $a_i = b_j = 0$, for
every $i$ and $j$. Then, we get $$ \sum_{i=1}^{12} c_i C_i +
\sum_{i=1}^6 d_i D_i + \sum_{i=1}^{12} e_i E_i + \sum_{i=1}^4
f_i F_i + \sum_{i=1}^3 g_i G_i = 0.  $$

\vspace{0.2cm}
\noindent
{\it Step 2.} Consider the homomorphism $\Fd\oticala P_4 \to
\Fd\oticala P_2$ induced by the projection $P_4 \to P_4/(z,
t) \cong P_2$. Under this homomorphism, the image of the above
linear relation is $d_1(5, 3) = 0$.

In order to show $d_1 = 0$, we need to prove that $(5, 3)$
is non zero in $\Fd\oticala P_2$.  The linear transformation
$x\mapsto x, y\mapsto x+y$ sents (5,3) to $(8, 0) + (7, 1) +
(6, 2) + (5, 3) \sim (7, 1) + (5, 3) $.  As the action of the
Steenrod algebra commutes with linear maps, if $(5, 3)$ is hit
then so is $(7, 1) + (5, 3)$. But it is impossible, because $(7,
1)$ is a spike. Hence, $(5, 3) \neq 0$ in $\Fd\oticala P_2$
and $d_1=0$.

Similarly, using all the projections of $P_4$ to its quotients
by the ideals generated by each pair of the four variables,
we get $d_i=0$ for every $i$. So we get $$ \sum_{i=1}^{12}
c_i C_i + \sum_{i=1}^{12} e_i E_i + \sum_{i=1}^4 f_i F_i +
\sum_{i=1}^3 g_i G_i = 0.  $$

\vspace{0.2cm}
\noindent
{\it Step 3.} Consider the homomorphism $\Fd\oticala P_4 \to
\Fd\oticala P_3$ induced by the projection $P_4 \to P_4/(t)
\cong P_3$. Under this homomorphism, the linear relation above
is sent to
$$
c_1(6, 1, 1) +c_4(1, 6, 1) + c_6(1, 1, 6) + e_1(5, 2, 1) + e_4(2, 5, 1) + e_6(2, 1, 5) = 0.
$$

Applying the linear map $x\mapsto x, y\mapsto x, z\mapsto y$ to
this relation, we obtain
\beqn
 (c_1+c_4+e_1+e_4) (7, 1) + c_6(2, 6) + e_6(3, 5) &=&\\
 (c_1+c_4+e_1+e_4) (7, 1) + e_6(3, 5) &=&  0.
\eeqn
Since $(7, 1)$ is a spike, $(c_1+c_4+e_1+e_4)=0$, hence $e_6(3,
5) = 0$. As for $(5, 3)$, we can show $(3, 5) \neq
0 \in \Fd\oticala P_2$ and get $e_6 =0$.

By similar arguments,
we have $e_1=e_4=e_6=0$.  The equality $(c_1+c_4+e_1+e_4)=0$
shows $c_1+c_4=0$ or $c_1=c_4$.  By similar arguments, $c_1=c_4=
c_6$. We denote this common coefficient by $c$ and get
$$
c \{
(6, 1, 1) +(1, 6, 1) + (1, 1, 6) \} = 0.
$$
We prove $c=0$ by showing $(6, 1, 1) +(1, 6, 1) + (1, 1, 6)
\neq 0$. Suppose the contrary, that $(6, 1, 1) +(1, 6, 1) +
(1, 1, 6)$ is hit. Then, by the unstable property of the action of
$\cala$ on the polynomial algebra, we have
$$
(6, 1, 1) +(1, 6, 1) + (1, 1, 6) = Sq^1(P) + Sq^2(Q) + Sq^4(R),
$$
for some
polynomials $P, Q, R$. By the degree information, $Sq^4(R)= R^2$
and this element is hit by $Sq^1$. Therefore, it suffices to
suppose $(6, 1, 1) +(1, 6, 1) + (1, 1, 6) = Sq^1(P) + Sq^2(Q)$.

\vspace{0.2cm}
Let $Sq^2Sq^2Sq^2$ act on the both sides of this equality. The right hand side is sent to zero, as
$Sq^2Sq^2Sq^2$ annihilates $Sq^1$ and $Sq^2$. On the other hand,
$$
Sq^2Sq^2Sq^2\{ (6, 1, 1) +(1, 6, 1) + (1, 1, 6) \} = (8, 4, 2) + \mbox{symmetries} \neq 0.
$$
This is a contradiction. So, it implies $(6, 1, 1) +(1, 6, 1) + (1, 1, 6) \neq 0$ and $c=0$. We get
$$
\sum_{i=1}^4 f_i F_i + \sum_{i=1}^3 g_i G_i = 0.
$$

\vspace{0.2cm}
\noindent
{\it Step 4.} Apply the linear map $x\mapsto x, y\mapsto y, z\mapsto y, t\mapsto y$ to the above
equality and we have
\beqn
   f_1(5, 3) + (f_2+f_3+f_4+g_3) (1, 7) + (g_1+g_2)(4, 4) &=&  \\
    f_1(5, 3) + (f_2+f_3+f_4+g_3) (1, 7)                 &=& 0.
\eeqn
As $(7, 1)$ is a spike, we obtain $(f_2+f_3+f_4+g_3)=0$ and $f_1(5, 3)=0$.
As $(5, 3) \neq 0$, it yields $f_1=0$.

Next, apply the linear map $x\mapsto x, y\mapsto y, z\mapsto x, t\mapsto x$ to the
equality $\sum_{i\neq 1} f_i F_i + \sum_{i=1}^3 g_i G_i = 0$
and we have
\beqn
  f_2(3, 5) + (f_3 + f_4 +g_2) (7, 1) + g_1(6, 2) + g_3(4, 4) &=& \\
  f_2(3, 5) + (f_3 + f_4 +g_2) (7, 1)  &=&  0.
\eeqn
As $(7, 1)$ is a spike, we get $(f_3 + f_4 +g_2)=0$ and $f_2(3, 5)= 0$.
Since $(3, 5)\neq 0$, it implies $f_2=0$.

Similarly, apply the linear map $x\mapsto x, y\mapsto x, z\mapsto y, t\mapsto x$ to the
equality $f_3F_3 + f_4F_4 + \sum_{i=1}^3 g_i G_i = 0$
and we have
\beqn
 f_3(3, 5) + (f_4 + g_1) (7, 1) + (g_2+g_3)(6, 2) &=& \\
 f_3(3, 5) + (f_4 + g_1) (7, 1)   &=&  0.
\eeqn
As $(7, 1)$ is a spike, we get $f_4+g_1=0$ and then $f_3=0$.

Finally, apply the linear map $x\mapsto x, y\mapsto x, z\mapsto x, t\mapsto y$ to the
equality $f_4F_4 + \sum_{i=1}^3 g_i G_i = 0$
and we have
\beqn
 f_4(3, 5) + (g_1 + g_2 + g_3)(7, 1) &=& \\
 f_4(3, 5) + (g_1 + g_2 + g_3)(7, 1) &=&  0.
\eeqn
As $(7, 1)$ is a spike, we get $g_1 + g_2 + g_3=0$ and then $f_4=0$.

Substituting $f_1=f_2=f_3=f_4=0$ into the equations $(f_2+f_3+f_4+g_3)=0$, $(f_3 + f_4 +g_2)=0$,
$f_4+g_1=0$, we get $g_1=g_2=g_3=0$.

We have shown that all coefficients of an arbitrary linear relation between the $55$ elements
listed in Proposition~\ref{basis} are zero.  The lemma follows.
\end{pf}

Combining Lemmas~\ref{generator} - \ref{independent}, we get Proposition~\ref{basis}.

\begin{Proposition}\label{0}
$(\Fd\oticala P_4)^{GL_4}_8=0.$
\end{Proposition}

\begin{pf}
If $X$ is one of the letters $A, B, C, D, E, F, G$, let $\call (X)$
be the vector subspace of
$(\Fd\oticala P_4)_8$ spanned by the elements of family $X$ in
Proposition~\ref{basis}.
Let $S_k$ denote the symmetric subgroup of $GL_k$.
According to the relations listed in the proof of Lemma~\ref{generator},
$\call (A)$, $\call (B)$, $\call (C)$, $\call (D)$, $\call (F)$, $\call (G)$ are $S_4$-submodules.
The subspace $\call (E)$ is not an $S_4$-submodule. However, the sum
$$
\call (C, E) = \call (C) \oplus \call (E)
$$
is.  We have a decomposition of $S_4$-modules
$$
(\Fd\oticala P_4)_8 = \call (A) \oplus \call (B) \oplus  \call (C, E) \oplus \call (D)
\oplus  \call (F) \oplus \call (G).
$$

Let $\alpha$ be an arbitrary $GL_4$-invariant in $(\Fd\oticala P_4)_8$. It can uniquely be written
in the form
$$
\alpha = \alpha_A + \alpha_B + \alpha_{C,E} + \alpha_D + \alpha_F + \alpha_G,
$$
where $\alpha_X \in \call(X)$ for $X \in \{ A, B, D, F, G\}$,
and $\alpha_{C, E} \in \call(C, E)$. Each term of this sum is $S_4$-invariant.

Note that if a linear combination of elements in a family
is $S_4$-invariant, then all of its coefficients are equal,
because each element in the family can be obtained from any
other by a suitable permutation. Let $s_X$ denote the sum
of all the elements in the family $X$ listed in
Proposition~\ref{basis}. Then, we have $\alpha_A =
a s_A$, $\alpha_B = b s_B$, $\alpha_D = d s_D$, $\alpha_F =
f s_F$, $\alpha_G = g s_G$, and $\alpha_{C, E} = c s_C +
e s_E$, where $a, b, c, d, e, f, g \in \Fd$.

Let $p$ be the transposition given by $p(x)=y$, $p(y)=x$, $p(z)=z$, $p(t)=t$. It is easy to see that
\beqn
p(2, 1, 0, 5) &=& (1, 2, 0, 5) = (2, 1, 0, 5) + (1, 1, 0, 6), \\
p(2, 1, 5, 0) &=& (1, 2, 5, 0) = (2, 1, 5, 0) + (1, 1, 6, 0).
\eeqn
Further, the 10 elements different from $(2, 1, 0, 5)$ and $(2, 1, 5, 0)$ in family $E$ are divided
into 5 pairs with $p$ acting on each pair by twisting.
So,  $p(s_E) = s_E + (1, 1, 0, 6) + (1, 1, 6, 0)$.
On the other hand, as the family $C$ is full, in the sense that it contains all the variable
permutations of a monomial, we have $p(s_C)= s_C$. Hence, we get
$$
p(\alpha_{C, E}) = p(c s_C + e s_E) = c s_C + e s_E + e(1, 1, 0, 6) + e(1, 1, 6, 0).
$$
As $\alpha_{C, E}$ is $S_4$-invariant, $e(1, 1, 0, 6) + e(1, 1, 6, 0)=0$. So $e=0$, because the two
elements are linearly independent by Lemma ~\ref{independent}. We obtain
$$
\alpha = \alpha_A + \alpha_B + \alpha_{C} + \alpha_D + \alpha_F + \alpha_G,
$$
where $\alpha_C= \alpha_{C,E}= c s_C$.

Let us now consider the transvection $\varphi$ given by
$\varphi (x)=x$, $\varphi (y)=y$, $\varphi (z)=z$, $\varphi (t)=x+t$. A routine computation shows
\beqn
\varphi(s_A) &=& s_A + (7, 1, 0, 0) + (7, 0, 1, 0)+( 7, 0, 0, 1 )+ (1, 7, 0, 0) + ( 1, 0, 7, 0 ) \\
             & & + ( 6, 1, 0, 1 ) + ( 6, 0, 1, 1 ) + ( 1, 1, 0, 6 ) + ( 1, 0, 1, 6 ),  \\
\varphi(s_B) &=& s_B +(6, 1, 1, 0)+(1, 6, 1, 0)+(1, 1, 6, 0)+(2, 5, 1, 0)+(2, 1, 5, 0)\\
             & & +(5, 1, 1, 1)+ (1, 5, 1, 1)+(1, 1, 5, 1) + (4, 2, 1, 1) + (4, 1, 2, 1) \\
             & & + (3, 3, 1, 1) + ( 3, 1, 3, 1 ), \\
\varphi(s_C) &=& s_C + ( 6, 1, 1, 0 ) + ( 1, 6, 1, 0 ) +( 1, 1, 6, 0 ), \\
\varphi(s_D) &=& s_D + (7, 0, 0, 1) + (1, 6, 0, 1) + (1, 0, 6, 1) + (5, 3, 0, 0) + (5, 0, 3, 0),\\
\varphi(s_F) &=& s_F + (2, 5, 1, 0) + (2, 1, 5, 0)+ (5, 1, 1, 1) + (4, 2, 1, 1) +( 4, 1, 2, 1), \\
\varphi(s_G) &=& s_G + ( 6, 1, 1, 0 ).
\eeqn

Let $r_X = \varphi(s_X) - s_X$ where $X$ is one of the letters
$A, B, C, D, F, G$. The equality
$\varphi(\alpha) = \alpha$ is rewritten as
$$
\varphi(a s_A + bs_B + cs_C + d s_D + f s_F + g s_G)=a s_A + bs_B + cs_C + d s_D + f s_F + g s_G,
$$
or equivalently
$$
a r_A + br_B + cr_C + d r_D + f r_F + g r_G =0.
$$

In this linear combination,  $r_B$ and $r_D$ are the only terms containing
$(3, 3, 1, 1)$ in family $B$ and $(5, 3, 0, 0)$ in family $D$ respectively.
From Lemma~\ref{independent},
we get $b= d = 0$, and therefore $a r_A + cr_C + f r_F + g r_G =0$.

%Similarly, as $r_F$ is the only term among $r_A, r_C, r_F, r_G$ which contains
% $(4, 2, 1, 1)$ in family $F$, it yields $f=0$, and then $a r_A + cr_C + g r_G =0$.

In the new linear combination, as $r_A$, $r_C$ and $r_F$ are the
only terms containing $(7, 1, 0, 0)$
in family $A$, $( 1, 6, 1, 0 )$ in family $C$ and  $(4, 2, 1, 1)$ in
family $F$ respectively, we have
$a= c  = f =0$.  As a consequence, $g r_G =0$, so we finally get $g=0$.

In summary, we have shown that every $GL_4$-invariant $\alpha$ in $(\Fd\oticala P_4)_8$ equals zero.
The proposition is proved.

\end{pf}

%------------------------------------
\section{The fourth algebraic transfer  is not an epimorphism}
         % does not detect the $g$-family elements}
\label{notepi}

The goal of this paper is to prove the following theorem, which is also numbered as
Therem~\ref{main}.

\begin{Theorem}
$$
Tr_4 : \Fd\otiglf PH_i(B\V_4) \to Ext_{\cala}^{4,4+i}(\Fd,\Fd)
$$
does not detect the non zero elements $g_s\in Ext_{\cala}^{4,12\cdot 2^s}(\Fd,\Fd)$
for every $s\geq 1$.
\end{Theorem}

\begin{pf} %{Proof of Theorem~\ref{main} }
Combining Lemma~\ref{8-string} and Proposition~\ref{0} we get
$$
(\Fd\oticala P_4)^{GL_4}_{12\cdot 2^s-4} = 0,
$$
for every non negative integer $s$.

On the other hand, it is well known that $Ext_{\cala}^{4, 24}(\Fd, \Fd)$ is spanned by
the generator $g_1$ (see May~\cite{May}). Further, $g_s = (Sq^0)^{s-1}(g_1)$ is non zero in
$Ext_{\cala}^{4,12\cdot 2^s}(\Fd,\Fd)$ (see Lin~\cite{Lin} and also ~\cite{LM}).

As $\Fd\otiglf PH_{12\cdot 2^s-4}(B\V_4)$ is dual
to $(\Fd\oticala P_4)^{GL_4}_{12\cdot 2^s-4}$,
$$
Tr_4 : \Fd\otiglf PH_{12\cdot 2^s-4}(B\V_4) \to Ext_{\cala}^{4,12\cdot 2^s}(\Fd,\Fd)
$$
does not detect the generator $g_s$, for every non negative integer $s$.

The theorem is proved.
\end{pf}

As a consequence, we get a negative answer to a prediction by Minami~\cite{Minami}.
(This corollary is also numbered as Corollary~\ref{localized}.)

\begin{Corollary}
The localization of the fourth algebraic transfer
$$
(Sq^0)^{-1}Tr_4 : (Sq^0)^{-1}\Fd\otiglf PH_i(B\V_4) \to (Sq^0)^{-1}Ext_{\cala}^{4,4+i}(\Fd,\Fd)
$$
given by inverting $Sq^0$ is not an epimorphism.
\end{Corollary}

\begin{pf}
Indeed, it does not detect the non zero element $g$, which is represented by the family
$(g_s)_{s>0}$ with $g_s = (Sq^0)^{s-1}(g_1)$. The corollary follows.
\end{pf}

\begin{Remark}
{\rm Our result does not affect Singer's conjecture
that the $k$-th algebraic transfer
is a monomorphism for every $k$. (See ~\cite{Singer89}.)
}
\end{Remark}

%-------------------------------------
%\begin{section}
\section{Final Remark: $GL_4$-module structure}

Boardman's study of the 3 variable problem shows that the
$GL_k$ module
structure of $\Fd\oticala P_k$
may be a useful tool.
In this vein we close with a description of the module
$(\Fd\oticala P_4)_8$  as a $GL_4$-module.  From the ``Modular Atlas''
\cite{ModAtlas} we find that there are 8 irreducible
modules for $GL_4$ in characteristic 2, of dimensions
1, 4, 4, 6, 14, 20, 20, and 64.
With a little
calculation we find the following description of them:
\begin{itemize}
\item[$~1$:] the trivial module $\Fd$
\item[$N$:] the natural module $\Fd^4$,
\item[$N^*$:] the dual of the natural module,
\item[$\Lambda$:] the alternating square of $N$ or $N^*$,
\item[$S$:] the nontrivial constituent of $N \otimes N^*$,
which has composition factors $1, S, 1$,
\item[$T$:] a constituent of $N \otimes \Lambda$,
which has composition factors $N^*$ and $T$,
\item[$T^*$:] a constituent of $N^* \otimes \Lambda$,
which has composition factors $N$ and $T^*$,
\item[$St$:] the Steinberg module.
\end{itemize}

Using a ``meataxe'' program written in MAGMA,
together with a MAGMA program to compute Brauer characters,
we have found
that $(\Fd \oticala P_4)_8$ is an extension
\[
0 \longrightarrow
N^* \oplus T \longrightarrow
(\Fd \oticala P_4)_8 \longrightarrow
\Lambda \oplus M \longrightarrow 0,
\]
where the 25-dimensional module $M$ is an extension
\[
0 \longrightarrow
1 \oplus \Lambda \longrightarrow
M \longrightarrow
N \oplus S \longrightarrow 0.
\]
The corresponding lattice
of submodules of $(\Fd\oticala P_4)_8$ is
shown in Figure~\ref{lattice}.
We name the submodules by their dimension,
using a prime to distinguish the two submodules
of dimension 30.
We label the edges 
by the corresponding quotient module.
In it, intersections are shown, but sums are omitted for
clarity.  That is, the intersection of the submodules
$30'$ and $35$ is the submodule $24$, but the sum of
$30'$ and $35$ (a submodule of dimension 41) is
not shown.
The two extensions above can be seen in the lattice,
in the sense that, for example, the submodule of dimension
24 is the direct sum of the submodules of dimensions 4 and 20, 
since their intersection is trivial.
Further, the quotient of 55 by 24 is the direct sum of
the quotients of 30' by 24 and of 49 by 24. 

\begin{figure}
\[
\xymatrix@!{
& 55 \ar@{-}[ldddd]_{M}   \ar@{-}[rrd]^{\Lambda}
& & &\\
& & & 49  \ar@{-}[ld]_{N} \ar@{-}[rd]^{S}
& \\
& & 45 \ar@{-}[rd]^{S}
& & 35 \ar@{-}[ld]_{N}
\\
& & & 31 \ar@{-}[rd]^{1}  \ar@{-}[ld]_{\Lambda}
& \\
30' \ar@{-}[rrrd]_{\Lambda}
& & 25 \ar@{-}[rd]^{1}
& & 30 \ar@{-}[ld]_{\Lambda}
 \\
& & & 24  \ar@{-}[rd]^{N^*} \ar@{-}[ld]_{T}
 & \\
& & 4 \ar@{-}[rd]^{N^*}
& & 20 \ar@{-}[ld]_{T}
 \\
& & & 0
 & \\
}
\]
\caption{Some $GL_4$-submodules of $(\Fd\oticala P_4)_8$.}
\label{lattice}
\end{figure}

The generators for these submodules are provided by the same computer
program used to find this decomposition and are listed below.    When
all the monomials in one of the seven families listed in
Proposition~\ref{basis} appear, we simply write the name of the family,
so that, for example, all the monomials in family $A$ are in  the
submodule of dimension 20.   Also, recall the element
\[
s_G =
( 4, 2, 1, 1 )+( 4, 1, 2, 1 )+( 1, 4, 2, 1 )
\]
used in the proof of Proposition~\ref{0}.  Finally note that
elements which form bases for the subquotients can
be read off by comparing these lists of generators. 
For example, the quotient of the module 30 by the submodule 24 is
$\Lambda$, and the elements of family $D$ generate it.

\begin{itemize}

\item[4:]

$( 6, 1, 1, 0 ) + ( 1, 6, 1, 0 ) + ( 1, 1, 6, 0 )$,
$( 6, 1, 0, 1 ) + ( 1, 6, 0, 1 ) + ( 1, 1, 0, 6 )$, \\
$( 6, 0, 1, 1 ) + ( 1, 0, 6, 1 ) + ( 1, 0, 1, 6 )$,
$( 0, 6, 1, 1 ) + ( 0, 1, 6, 1 ) + ( 0, 1, 1, 6 )$,

\item[20:]  $(A)$,
$( 6, 1, 1, 0 ) + ( 1, 1, 6, 0 )$,
$( 6, 1, 0, 1 ) + ( 1, 1, 0, 6 )$,
$( 6, 0, 1, 1 ) + ( 1, 0, 1, 6 )$,  \\
$( 1, 6, 1, 0 ) + ( 1, 1, 6, 0 )$,
$( 1, 6, 0, 1 ) + ( 1, 1, 0, 6 )$,
$( 1, 0, 6, 1 ) + ( 1, 0, 1, 6 )$, \\
$( 0, 6, 1, 1 ) + ( 0, 1, 1, 6 )$,
$( 0, 1, 6, 1 ) + ( 0, 1, 1, 6 )$.

\item[24:] $(A)$ and $(C)$.

\item[25:]
 $(A),  (C)$, and $s_G$.

\item[30:]
 $(A),  (C)$, and $(D)$.

\item[$30'$:] $(A)$, $(C)$ and
$( 5, 1, 1, 1 ) + ( 1, 5, 1, 1 ) + s_G +  ( 3, 3, 1, 1 )  $,\\
$( 5, 1, 1, 1 ) + ( 1, 1, 5, 1 ) + s_G + ( 3, 1, 3, 1 )  $,
$( 5, 1, 1, 1 ) + ( 1, 1, 1, 5 ) + s_G + ( 3, 1, 1, 3 )  $,\\
$( 1, 5, 1, 1 ) + ( 1, 1, 5, 1 ) + s_G + ( 1, 3, 3, 1 )  $,
$( 1, 5, 1, 1 ) + ( 1, 1, 1, 5 ) + s_G + ( 1, 3, 1, 3 )  $,\\
$( 1, 1, 5, 1 ) + ( 1, 1, 1, 5 ) + s_G + ( 1, 1, 3, 3 )  $.

\item[31:]
 $(A),  (C), (D)$ and $s_G$.

\item[35:]
 $(A),  (C), (D)$, $s_G$ and\\
$( 5, 2, 1, 0 )+( 5, 2, 0, 1 )+ ( 5, 0, 2, 1 ) + ( 5, 1, 1, 1 )$,\\
$( 2, 5, 1, 0 )+( 2, 5, 0, 1 )+ ( 0, 5, 2, 1 ) + ( 1, 5, 1, 1 )$,\\
$( 2, 1, 5, 0 )+( 2, 0, 5, 1 )+ ( 0, 2, 5, 1 ) + ( 1, 1, 5, 1 )$,\\
$( 2, 1, 0, 5 )+( 2, 0, 1, 5 )+ ( 0, 2, 1, 5 ) + ( 1, 1, 1, 5 )$.

\item[45:]
 $(A),  (C), (D), (E)$  and $(G)$.

\item[49:]
 $(A),  (C), (D), (E), (F)$  and $(G)$.

\end{itemize}
%\end{section}

%-------------------------------------

\newpage
\vspace{0.5cm}
\noindent
DEPARTMENT OF MATHEMATICS,  WAYNE STATE UNIVERSITY \newline
656 W. KIRBY STREET, DETROIT, MI 48202 (USA) \newline
ROBERT R. BRUNER: rrb@@math.wayne.edu  \newline
NGUY\EEX N H. V. H\UW NG: nhvhung@@math.wayne.edu

\vspace{0.5cm}
\noindent
IHES, F-91440, BURE-SUR-YVETTE, FRANCE  \newline
L\^{E} M. H\`{A}: lha@@ihes.fr  \newline
\end{document}